\documentclass[a4paper]{article}
\pdfoutput=1

\usepackage{graphicx,wrapfig}
\usepackage{lipsum}
\usepackage{amsmath,amssymb,mathrsfs}
\usepackage{hyperref}
\usepackage{url}
\usepackage{mathtools}
\usepackage{changepage}
\usepackage{multirow}
\usepackage{wrapfig}
\usepackage{subfig}
\usepackage{color}
\usepackage{epsfig,enumerate,amsmath,amsfonts,amssymb,amsthm,mathrsfs,ifpdf}
\usepackage{indentfirst,relsize}
\usepackage{setspace,graphicx}
\usepackage{latexsym}
\usepackage[all]{xy}
\usepackage[usenames,dvipsnames]{pstricks}
\usepackage{pst-grad} 
\usepackage{pst-plot} 
\usepackage[margin = 2.50cm]{geometry}

\usepackage{algorithm}

\usepackage[noend]{algpseudocode}

\newcommand{\remove}[1]{}

\newtheorem{theo}{Theorem}
\newtheorem{lem}[theo]{Lemma}

\newtheorem{cl}[theo]{Claim}

\graphicspath{{./images/}}

\title{On $k$-vertex-edge domination of graph}

\author{Debojyoti Bhattacharya\footnote{Inian Instittute of Technology Patna, Bihta, 801106, Bihar, India. email: debojyoti\_2021ma11@iitp.ac.in} \and Subhabrata Paul\footnote{Inian Instittute of Technology Patna, Bihta, 801106, Bihar, India. email: subhabrata@iitp.ac.in}}

\date{}

\begin{document}
	\maketitle
	
\begin{abstract}
	Let $G=(V,E)$ be a simple undirected graph. The open neighbourhood of a vertex $v$ in $G$ is defined as $N_G(v)=\{u\in V~|~ uv\in E\}$; whereas the closed neighbourhood is defined as $N_G[v]= N_G(v)\cup \{v\}$. For an integer $k$, a subset $D\subseteq V$ is called a $k$-vertex-edge dominating set of $G$ if for every edge $uv\in E$, $|(N_G[u]\cup N_G[v]) \cap D|\geq k$. In $k$-vertex-edge domination problem, our goal is to find a $k$-vertex-edge dominating set of minimum cardinality of an input graph $G$. In this paper, we first prove that the decision version of $k$-vertex-edge domination problem is NP-complete for chordal graphs. On the positive side, we design a linear time algorithm for finding a minimum $k$-vertex-edge dominating set of tree. We also prove that there is a $O(\log(\Delta(G)))$-approximation algorithm for this problem in general graph $G$, where $\Delta(G)$ is the maximum degree of $G$. Then we show that for a graph $G$ with $n$ vertices, this problem cannot be approximated within a factor of $(1-\epsilon) \ln n$ for any $\epsilon >0$ unless $NP\subseteq DTIME(|V|^{O(\log\log|V|)})$. Finally, we prove that it is APX-complete for graphs with bounded degree $k+3$.

\noindent\textbf{keywords:} {$k$-vertex-edge domination, NP-completeness, Approximation Algorithm, APX-completeness}
\end{abstract}


\section{Introduction}\label{sec1}

%
Domination and its variations are considered to be one of the classical problems in graph theory due to its application in different areas. Let $G=(V,E)$ be a simple undirected graph. The open neighbourhood of a vertex $v$ in $G$ is defined as $N_G(v)=\{u\in V~|~ uv\in E\}$; whereas the closed neighbourhood is defined as $N_G[v]= N_G(v)\cup \{v\}$. A subset $D\subseteq V$ is called a dominating set of $G$ if for every $v\in V$, $|N_G[v]\cap D|\geq 1$. In general, a vertex dominates its neighbouring vertices and also itself. Depending on the nature of dominating power of a vertex, different variations of this classical domination problem have been studied in literature \cite{domingraphbook,Fundamentalbook}.

In one of the variations, a vertex $v$ is considered to dominate all the edges that are incident to any vertex in $N_G[v]$. In literature, it is referred to as \emph{vertex-edge domination} or \emph{ve-domination} in short. An edge $e=uv\in E$ is said to be vertex-edge dominated(\emph{ve-dominated}) by a vertex $x$ if $x\in N_G[u]\cup N_G[v]$. A subset $D_{ve}\subseteq V$ is called a \emph{vertex-edge dominating set} or \emph{ve-dominating set} of $G$ if for every edge $uv\in E$, $|(N_G[u]\cup N_G[v]) \cap D_{ve}|\geq 1$, that is, every edge of the graph is ve-dominated by $D_{ve}$. The minimum cardinality of a ve-dominating set of a graph $G$ is called \emph{ve-domination number} of $G$ and it is denoted as $\gamma_{ve}(G)$. In this problem, the goal is to find a ve-dominating set of minimum cardinality in a given input graph. The notion of a ve-dominating set was introduced by Peters in his Ph.D. thesis \cite{peters}. This problem has been well studied both from algorithmic as well as theoretical point of view \cite{boutrig2016vertex, chitra2012global, jena, krishnakumari2014bounds, lewis, paul2, paul, peters, naresh, zylinski2019vertex, chen2022double}. 

A generalization of vertex-edge domination, namely $k$-vertex-edge domination was studied by Li and Wang in 2023 \cite{li2023polynomial}. Given an integer $k$, a subset $D_{kve}\subseteq V$ is called a \emph{$k$-vertex-edge dominating set} or \emph{$k$-ve dominating set} of $G$ if for every edge $uv\in E$, $|(N_G[u]\cup N_G[v]) \cap D_{kve}|\geq k$. The minimum cardinality of a $k$-ve dominating set of a graph $G$ is called \emph{$k$-ve domination number} of $G$ and it is denoted as $\gamma_{kve}(G)$. The minimum $k$-ve domination problem and its corresponding decision version are defined as follows:

\noindent\underline{\textbf{Minimum $k$-Vertex-Edge Domination Problem}(\textsc{Min$k$VEDP})}

\noindent\emph{Instance}: A graph $G=(V,E)$ and an integer $k$.

\noindent\emph{Output}: Minimum $k$-vertex-edge dominating set $D$ of $G$.

\noindent\underline{\textbf{Decision version of $k$-Vertex-Edge Domination Problem}(\textsc{Decide$k$VEDP})}

\noindent\emph{Instance}: A graph $G=(V,E)$ and an integer $k$ and an integer $t$.

\noindent\emph{Question}: Does there exists a $k$-vertex-edge dominating set $D$ of $G$ of size at most $t$?

In \cite{li2023polynomial}, authors proposed a $O(km)$ time algorithm for \textsc{Min$k$VEDP} in interval graph and a linear time algorithm to find a minimum independent vertex-edge dominating set in unit interval graph. A set $S\subset V$ is an independent vertex-edge dominating set if it is independent and vertex-edge dominating set.    
The results presented in this paper are as follows:

\begin{itemize}
	\item We show that \textsc{Decide$k$VEDP} is NP-complete for chordal graphs.
	
	\item The \textsc{Min$k$VEDP} can be solved in linear time for trees.
	\item  A lower bound on approximation ratio of \textsc{Min$k$VEDP}.
	\item An $O(\log \Delta(G))$-approximation algorithm for \textsc{Min$k$VEDP}.
	\item The \textsc{Min$k$VEDP} is APX-hard for graphs with maximum degree $k+3$. 
\end{itemize}

The rest of the paper is organized as follows. In Section \ref{sec2}, we present the NP-completeness result for \textsc{Decide$k$VEDP} in chordal graphs. Section \ref{sec3} deals with the linear time algorithm for \textsc{Min$k$VEDP} in tree. After that, in Section \ref{sec4} we discuss approximation algorithms and hardness of approximation results for \textsc{Min$k$VEDP}. Finally, Section \ref{sec5} concludes this paper with some interesting open problems.

%
\section{NP-complete for chordal graphs}\label{sec2}

A graph $G=(V,E)$ is \emph{chordal} graph if every cycle of length greater or equal to $4$ has a chord. A chord is an edge between two non-consecutive vertices of a cycle. In this section, we show that \textsc{Decide$k$VEDP} is NP-complete for chordal graphs. We prove this result by reducing \textsc{Exact $3$ Cover Problem} (\textsc{Ex3CP}) to our problem. The \textsc{Ex3CP} is known to be NP-complete \cite{garey1979} and the problem is as follows:

\noindent\underline{\textsc{Exact $3$ Cover Problem} (\textsc{Ex3CP})}

\noindent\emph{Instance}: A set $X$ of $3q$ elements and a collection of $3$ element subsets of $X$, say $C$.

\noindent\emph{Question}: Does there exist a sub-collection $C'\subseteq C$ such that every element of $X$ belong to exactly one member of $C'$?

\begin{theo}
	The \textsc{Decide$k$VEDP} is NP-complete for chordal graphs.
\end{theo}
\begin{proof}
	Given a subset of $V$ of size $t$, we can verify whether it is a $k$-ve dominating set of $G$ or not in polynomial time. Therefore, the \textsc{Decide$k$VEDP} is in NP. Next, we show a polynomial time reduction from an instance of \textsc{Ex3CP} to an instance of \textsc{Decide$k$VEDP}. Let $(X, C)$ be an instance of \textsc{Ex3CP} where $X=\{x_1, x_2, \ldots , x_{3q}\}$ and $C=\{c_{1},c_{2},...,c_{m}\}$. We construct the graph $G$ as follows:
	\begin{figure}
		\centering
		\includegraphics[scale=0.6]{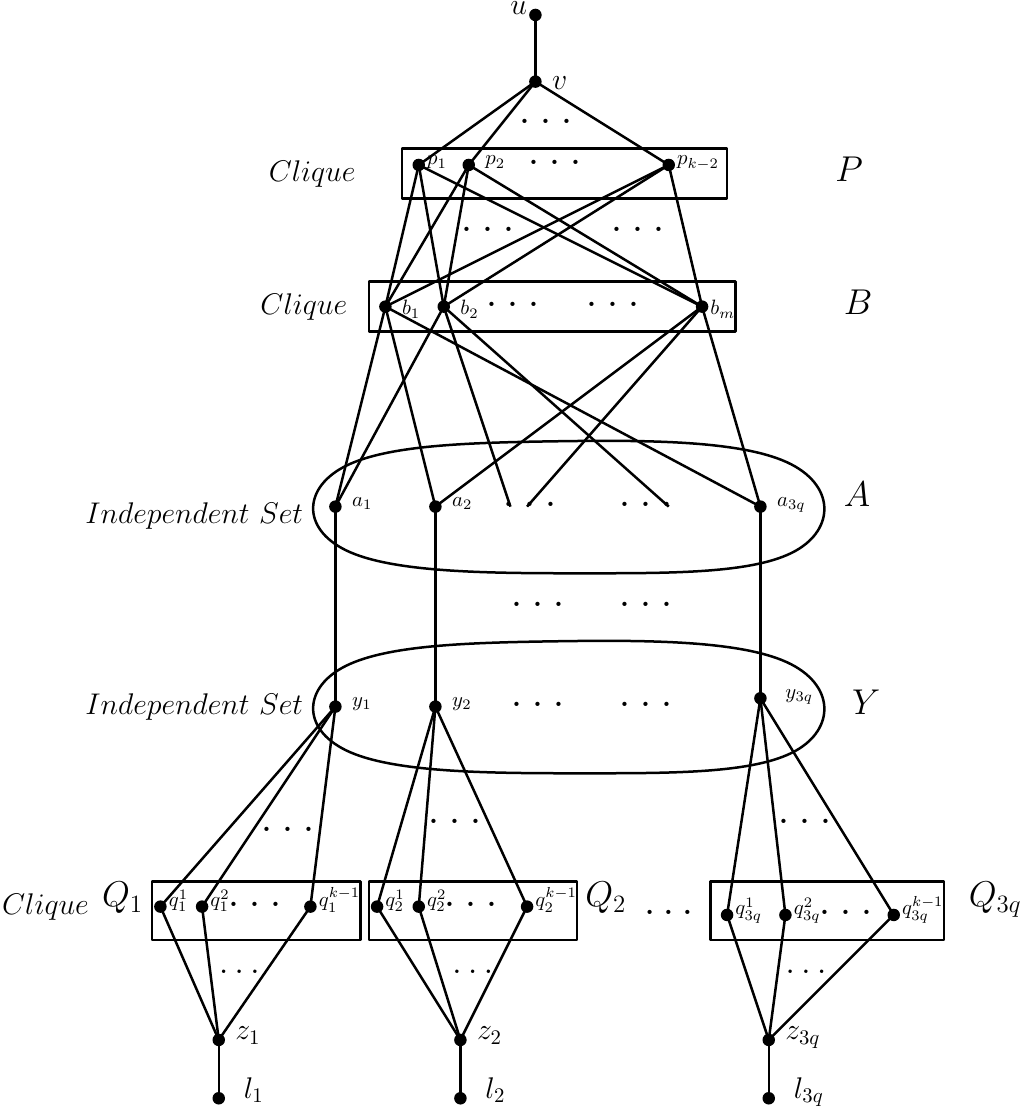}
		\caption{Construction of $G$ from an instance of \textsc{Ex3CP}}\label{Fig: NPC for chordal}
	\end{figure}
	
	\begin{itemize}
		\item For every $x_{i}\in X$, consider a vertex $a_{i}$ and for every $c_{j}\in C$, consider a vertex $b_{j}$ in $G$. Let $A=\{a_1, a_2, \ldots , a_{3q}\}$ and $B=\{b_{1},b_{2},\ldots,b_{m}\}$. We make $B$ a clique by adding all possible edges in $B$. Also, if $x_i\in c_j$, then we add an edge between $a_i$ and $b_j$.
		
		\item Consider a clique $P=\{p_1, p_2, \ldots p_{k-2}\}$ of size $k-2$ and make every vertex of $P$ adjacent to every vertex of $B$. Note that, $B\cup P$ is also a clique.
		
		\item Consider a path of length $2$, that is, $P_2=uv$ and add edge between $v$ and every vertex of $P$.
		
		\item For every vertex of $A$, consider a vertex $y_i$ and add the edge between $a_i$ and $y_i$ for every $i\in \{1, 2, \ldots, 3q\}$. Let $Y=\{y_1, y_2, \ldots , y_{3q}\}$.
		
		\item For every $i\in \{1, 2, \ldots, 3q\}$, consider a clique $Q_i$ of size $(k-1)$ and add edge between $y_i$ and every vertex of $Q_i$. 
		
		\item For every $i\in \{1, 2, \ldots, 3q\}$, consider a path of length $2$, that is, $P^i_2=z_il_i$ and add edge between $z_i$ and every vertex of $Q_i$.
	\end{itemize}
	Note that, in $G$, both $A$ and $Y$ form independent sets. Moreover, since each of $P$, $B$ and $Q_i$ for all $i\in \{1, 2, \ldots, 3q\}$ forms clique in $G$, there is no induced cycle of length more than $3$ in $G$. Hence, $G$ is a chordal graph. The construction of $G$ from the instance of \textsc{Ex3CP} is illustrated in Figure \ref{Fig: NPC for chordal}. To show the NP-completeness of the decision version of \textsc{Min-$k$-vedp} for chordal graphs, next we prove the following claim.
	
	\begin{cl}\label{Claim:NPC chordal}
		The instance $(X,C)$ of \textsc{Ex3CP} has an exact cover of size $q$ if and only if $G$ has a $k$-ve dominating set of size at most $k+q+3qk$.
	\end{cl}
	
	\begin{proof}[Proof of Claim \ref{Claim:NPC chordal}]
		Let $C'$ be an exact cover of $X$ of size $q$ and $B'=\{b_j~|~c_j\in C'\}$. Consider $D= \{u, v\} \cup P\cup B^{'}\bigcup\limits_{i=1}^{3q}(Q_{i}\cup\{z_{i}\})$. Clearly $|D|= k+q+3qk$. Note that the vertex set $(\{u, v\} \cup P\cup B^{'})$ $k$-ve dominates $uv$, $vp_i$ for every $p_i\in P$, every edge inside the clique $P\cup B$ and every edge of the form $a_ib_j$, since $B'$ has size at least $2$. Also, the vertex set $B^{'}\cup \bigcup\limits_{i=1}^{3q}Q_{i}$ $k$-ve dominates every edge of the form $a_iy_i$. Furthermore, the vertex set $\bigcup\limits_{i=1}^{3q}(Q_{i}\cup\{z_{i}\})$ $k$-ve dominates $z_il_i$, every edge between $z_i$ and $Q_i$, every edge inside $Q_i$ and every edge between $y_i$ and $Q_i$ for every $i\in \{1, 2, \ldots, 3q\}$. Hence $D$ is a $k$-ve dominating set of size $k+q+3qk$.
		
		Conversely, let $D$ be a $k$-ve dominating set of $G$ of size at most $k+q+3qk$. Note that to $k$-ve dominate the edge $uv$ , $\{u, v\}\cup P$ must be in $D$. On the other hand, for each $i \in \{1, 2, \ldots, 3q\}$, to dominate $z_il_i$ at least $k$ vertices from $(Q_i\cup \{z_i, l_i\})$ must be in $D$. In case where $D$ contains all $k+1$ vertices of $(Q_i\cup \{z_i, l_i\})$ for some $i$, we can remove $l_i$ from $D$ to obtain another $k$-ve dominating set of $G$ of size at most $k+q+3qk$. Therefore, without loss of generality, let us assume that $D$ contains  $(Q_i\cup \{z_i\})$ for every $i \in \{1, 2, \ldots, 3q\}$. Also, note that $D$ must contain vertices from $B$. Because every edge of the form $a_iy_i$ is $(k-1)$-ve dominated by $(\{u, v\} \cup P \bigcup\limits_{i=1}^{3q}(Q_{i}\cup\{z_{i}\}))$ and to $k$-ve dominate $a_iy_i$, we need at least one vertex from $\{a_i, y_i\}\cup B$. If $D$ does not contain any vertex from $B$, then the size of $D$ would be more than $k+q+3qk$ because we require at least $3q$ vertices from $A\cup Y$. Therefore, $D$ contains vertices from $B$. Also, observe that if $D$ contains both $a_i$ and $y_i$ for some $i\in \{1, 2, \ldots, 3q\}$, then we can remove $y_i$ from $D$ to get another $k$-ve dominating set of $G$ of size at most $k+q+3qk$. Moreover, if for some $i\in \{1, 2, \ldots, 3q\}$, $a_i\notin D$ but $y_i\in D$, then we can replace $y_i$ by $a_i$ to get another $k$-ve dominating set of $G$ of size at most $k+q+3qk$. Therefore, without loss of generality, let us assume that $D$ does not contain any vertex from $Y$. Let $D'= D \setminus (\{u, v\} \cup P \bigcup\limits_{i=1}^{3q}(Q_{i}\cup\{z_{i}\}))$. Clearly $|D'|\leq q$. Also, note that $D'$ contains only vertices from $A \cup B$. Let us assume that $D'$ contains $t_A$ and $t_B$ many vertices from $A$ and $B$, respectively. Now, note that $t_A$ many vertices from $A$ can dominate at most $t_A$ many edges of the form $a_iy_i$ and $t_B$ many vertices from $B$ can ve-dominate at most $3 t_B$ many edges of the form $a_iy_i$. Therefore, $t_A + 3 t_B \geq 3q$. Also, we know that $t_A + t_B \leq q$. This implies that $t_A=0$ and $t_B = q$. Therefore $D'\subset B$. Let $C' = \{c_j~|~b_j\in D'\}$ be a sub-collection of $C$. Clearly, $|C'|=q$. Since, $D'$ ve-dominates every edge of the form $a_iy_i$ and $|C'|=q$, every element of $X$ belongs to exactly one member of $C'$. Therefore, $C'$ is a solution of \textsc{Ex3CP}.
	\end{proof}
	From the above claim, it follows that the \textsc{Decide$k$VEDP} is NP-complete for chordal graphs.
\end{proof}

\section{Algorithm for trees}\label{sec3}
In this section, we give a linear time algorithm for the \textsc{Min$k$VEDP} in tree. Our algorithm is a greedy algorithm and based on the tree ordering and labeling method described in \cite{labelingC}.
\subsection{Algorithm}
Let $T=(V,E)$ be a tree. For every $v\in V $, we define the following label- $t(v)=\{B,R\}$. Thus, we partition the vertex set into two sets $B$ and $R$. Also, we assign some nonnegative integers $s(e)$ to every edge $e\in E $. Instead of determining the minimum $k$-ve dominating set, our algorithm determines a minimum $(s,t)$-dominating set.

The $(s,t)$-dominating set $D$ is defined as follows-\\
$(1)$ If $t(v)=R$, then $v\in D$.\\
$(2)$ For every edge $e=uv$, $|(N_{T}[u]\cup N_{T}[v])\cap D|\geq s(e)$.

If $t(v)=B$ for every $v\in V$ and $s(e)=k$ for every $e\in E$, then the $(s,t)$-ve dominating set is the $k$-ve dominating set. Let $v$ be a leaf of $T$, $u$ be the parent of $v$, $w$ be the parent of $u$, $r(u)$ be the set of vertices in $N_{T}[u]$ labeled as $R$ and $c(u)$ is the set of leaves of $u$.

\begin{algorithm}[H]
	\scriptsize
	\caption{\small VEDS-TREE$(T,s,t)$}
	\label{Algo:k-VEDS_tree}
	\textbf{Input:} $T=(V,E)$ and integers $s(e) ~for ~every ~ e\in E$ and $t(v)=\{B,R\}~for~every~ v\in V$.\\
	\textbf{Output:} A minimum $(s,t)$-dominating set $D$.
	\begin{algorithmic}[1]
		\State Find an ordering $\sigma=\{v_1,v_2,\ldots,v_n\}$ such that $\{v_n,v_{n-1},\ldots,v_1\}$ is the BFS ordering;
		\State $D=\phi$;
		\For {$(every ~support~vertex ~u\in \sigma)$}
		\State Find a vertex $v\in c(u)$ such $s(uv)=\max\{s(uz)|z\in N_{T}(u)\setminus\{w\}\}$;
		\State Compute $r(u)$;
		\If{$(s(uv)>|N_{T}[u]|~or~s(uw)>|N_{T}[u]\cup N_{T}[w]|)$}
		\State STOP.
		\ElsIf{$(s(uv) = |N_{T}[u]|$ or $s(uw)=|N_{T}[u]\cup N_{T}[w]|)$}
		\State Relabel $t(w)=R$ and $t(u)=R$; 
		\State $s(uw)=\max\{s(uw)-|c(u)|,0\}$;
		\State $D=D\cup c(u)$;
		\Else
		\If{$(s(uv)\leq |r(u)|)$}
		\State $s(uw)=\max\{s(uw)-|c(u)\cap r(u)|,0\}$;
		\State $D=D\cup (c(u)\cap r(u))$;
		\Else
		\If{$((s(uv)-|c(u)\cap r(u)|)=1)$}
		\State Relabel $t(w)=R$;
		\State $s(uw)=\max\{(s(uw)-|c(u)\cap r(u)|),0\}$;
		\State $D=D\cup (c(u)\cap r(u))$;
		\ElsIf{$(t(u)=B)$}
		\State Relabel $t(w)=R$ and $t(u)=R$;
		\State $s(uw)=\max\{(s(uw)-s(uv)+2),0\}$;
		\State $D=D\cup(c(u)\cap r(u))$;
		\State Include 
		$(s(uv)-|c(u)\cap r(u)|-2)$ many vertices of $c(u)\setminus r(u)$ in $D$;
		\Else
		\State Relabel $t(w)=R$;
		\State $s(uw)=\max\{(s(uw)-s(uv)+2),0\}$;
		\State $D=D\cup (c(u)\cap r(u))$;
		\State Include $(s(uv)-|c(u)\cap r(u)|-2)$ many vertices of $c(u)\setminus r(u)$ in $D$;
		\EndIf					
		\EndIf
		\EndIf
		\State $T=T\setminus c(u)$;
		\EndFor
		\State \Return $D$;
	\end{algorithmic}
\end{algorithm} 
We first apply BFS to the tree $T$ to find an ordering $\sigma =\{v_1,v_2,\ldots ,v_n\}$ where $v_n$ is the root of the tree. Starting with a support vertex $u$ of $T$ in $\sigma$, we proceed by checking the $s$-values of the incident edges of $u$ and find the pendent edge with maximum $s(e)$ value and corresponding leaf $v\in c(u)$. 
Depending on $s(uv)$ and the label of the vertices in $N_{T}[u]$, we either include $c(u)$ in $D$ and then delete $c(u)$ from $T$ or include some vertices of $c(u)$ in $D$, relabel the vertices $(w,u)$ and delete $c(u)$ from $T$. Also, we update the $s$-value of the edge $uw$. 
Details are in the above algorithm. 
\subsection{Proof of correctness}
Now, we show the correctness of our algorithm. Let $v$ be a leaf of $T$, $u$ be the parent of $v$ and $w$ be the parent of $u$. Let $\gamma_{(s,t)}(T)$ denote the $(s,t)$-ve domination number of $T$. Let $v\in c(u)$ such that $s(uv)=\max\{s(uz)|z\in N_{T}(u)\setminus\{w\}\}$.
\begin{lem}
	If $s(uv) > |N_{T}[u]|$ or $s(uw)>|N_{T}[u]\cup N_{T}[w]|$ then there is no dominating set.
\end{lem}
\begin{proof}
	Follows directly from the definition.
\end{proof}

\begin{lem}
	If $s(uv) = |N_{T}[u]|$ or $s(uw)=|N_{T}[u]\cup N_{T}[w]|$, then $\gamma_{(s,t)}(T)=\gamma_{(s',t')}(T')+|c(u)|$, where $T'$ is obtained from $T$ by deleting $c(u)$ and by relabelling $t'(u)=t'(w)=R$ and $s'(uw)= \max\{(s(uw)-|c(u)|),0\}$ and every other label remains the same.
\end{lem}
\begin{proof}
	Let $D$ be a minimum $(s,t)$-dominating set of $T$.	If $s(uv) = |N_{T}[u]|$ or $s(uw)=|N_{T}[u]\cup N_{T}[w]|$, then $D$ contains every vertex of $N_{T}[u]$. Let $D'=D\setminus c(u)$. We show that $D'$ is an $(s',t')$-dominating set of $T'$. Since $D$ contains every vertex of $N_{T}[u]$, the vertices $u$ and $w$ is in $D'$. For the edge $uw$, we know that $|(N_{T}[u]\cup N_{T}[w])\cap D|\geq s(uw)$ in $T$. This implies that $|(N_{T'}[u]\cup N_{T'}[w])\cap D'|\geq \max\{(s(uw)-|c(u)|),0\} = s'(uw)$ in $T'$. Since every other label remains the same, for every other edge $xy$, $|(N_{T'}[x]\cup N_{T'}[y])\cap D'|\geq s'(xy)$. Hence, $D'$ is an $(s',t')$-dominating set of $T'$. Therefore, $\gamma_{(s',t')}(T')\leq \gamma_{(s,t)}(T)-|c(u)|$.
	
	On the other hand, let $D'$ be a minimum $(s',t')$-dominating set of $T'$. Since $t'(u)=t'(w)=R$, $u,w\in D'$. Let $D= D' \cup c(u)$. For the edge $uw$, we know that $|(N_{T'}[u]\cup N_{T'}[w])\cap D'| \geq s'(uw)$. This implies that $|(N_{T}[u]\cup N_{T}[w])\cap D| \geq s'(uw)+ |c(u)|\geq s(uw)$. Also, for the other edge $uz$ incident on $u$ in $T$, we have $|(N_{T}[u]\cup N_{T}[z])\cap D|\geq |N_{T}[u]| \geq s(uz)$. Since every other label remains the same, for every other edge $xy$, $|(N_{T}[x]\cup N_{T}[y])\cap D|\geq s(xy)$. Hence, $D$ is an $(s,t)$-dominating set of $T$ and $\gamma_{(s,t)}(T)\leq \gamma_{(s',t')}(T')+|c(u)|$. Therefore we have,
	$\gamma_{(s,t)}(T)= \gamma_{(s',t')}(T')+|c(u)|$.
\end{proof}

\begin{lem}
	Let $s(uv) < |N_{T}[u]|$. If $s(uv)\leq |r(u)|$, then $\gamma_{(s,t)}(T)=\gamma_{(s',t')}(T')+|c(u)\cap r(u)|$, where $T'$ is obtained from $T$ by deleting $c(u)$ and by relabelling $s'(uw)= \max\{(s(uw)-|c(u)\cap r(u)|),0\}$ and every other label remains the same.
\end{lem}
\begin{proof}
	Let $D$ be a minimum $(s,t)$-dominating set of $T$.	Therefore, $r(u)\subseteq D$. Let $D'=D\setminus (c(u)\cap r(u))$. Now, we show that $D'$ is an $(s',t')$-dominating set of $T'$. For the edge $uw$, we have $|(N_{T}[u]\cup N_{T}[w])\cap D|\geq s(uw)$. This implies that $|(N_{T'}[u]\cup N_{T'}[w])\cap D'|\geq \max\{(s(uw)-|c(u)\cap r(u)|,0)\}=s'(uw)$. Since every other label remains the same, for every other edge $xy$, $|(N_{T'}[x]\cup N_{T'}[y])\cap D'|\geq s'(xy)$. Hence, $D'$ is an $(s',t')$-dominating set of $T'$. Therefore, $ \gamma_{(s',t')}(T')\leq\gamma_{(s,t)}(T)-|c(u)\cap r(u)|$.
	
	Also, let $D'$ be a minimum $(s',t')$-dominating set of $T'$. Let $D=D'\cup (c(u)\cap r(u))$. For the edge $uw$, we know that $|(N_{T'}[u]\cup N_{T'}[w])\cap D'|\geq s'(uw)$. This implies that $|(N_{T}[u]\cup N_{T}[w])\cap D|\geq s'(uw)+|c(u)\cap r(u)| \geq s(uw)$.
	Observe that, $D'$ contains every vertex of $N_{T'}[u]$ with label $R$.
	Since $s(uv)\leq |r(u)|$, for every other edge $uz$ incident on $u$ in $T$, we have $|(N_{T}[u]\cup N_{T}[z])\cap D|\geq |r(u)| \geq s(uz)$.
	Since every other label remains the same, for every other edge $xy$, $|(N_{T}[x]\cup N_{T}[y])\cap D|\geq s(xy)$. Thus, $D$ is an $(s,t)$-dominating set in $T$. Hence, $\gamma_{(s,t)}(T)\leq\gamma_{(s',t')}(T')+|c(u)\cap r(u)|$. Therefore, we have $\gamma_{(s,t)}(T)=\gamma_{(s',t')}(T')+|c(u)\cap r(u)|$.
\end{proof}

\begin{lem}
	Let $s(uv) < |N_{T}[u]|$ and $s(uv)>|r(u)|$. If $(s(uv)- |c(u)\cap r(u)|)=1$, then $\gamma_{(s,t)}(T)=\gamma_{(s',t')}(T')+|c(u)\cap r(u)|$, where $T'$ is obtained from $T$ by deleting $c(u)$ and relabelling $t'(w)=R$ and $s'(uw)= \max\{(s(uw)-|c(u)\cap r(u)|),0\}$ and every other label remains the same.
\end{lem}

\begin{proof}
	Let $D$ be a minimum $(s,t)$-dominating set of $T$. Therefore, $r(u)\subseteq D$. Since $(s(uv)- |c(u)\cap r(u)|)=1$ and $s(uv)>|r(u)|$, we have $t(w)=B$. Moreover, since $D$ is an $(s,t)$-dominating set of $T$,  there must be a vertex $z\in D$ such that $z\in N_{T}[u]$ with $t(z)=B$. Let $D'=(D\setminus\{z\})\cup\{w\}$.
	Clearly, $D'$ is also a minimum $(s,t)$-dominating set of $T$. Let $D''=D'\setminus(c(u)\cap r(u))$. We show that $D''$ is an $(s',t')$-dominating set of $T'$. Clearly, $w\in D''$. For the edge $uw$, we know that $|(N_{T}[u]\cup N_{T}[w])\cap D'|\geq s(uw)$. This implies that $|(N_{T'}[u]\cup N_{T'}[w])\cap D''|\geq \max\{(s(uw)-|c(u)\cap r(u)|),0\}=s'(uw)$. Since every other label remains the same, for every other edge $xy$, $|(N_{T'}[x]\cup N_{T'}[y])\cap D''|\geq s'(xy)$. Hence, $D''$ is an $(s',t')$-dominating set of $T'$. Therefore, $\gamma_{(s',t')}(T')\leq\gamma_{(s,t)}(T)-|c(u)\cap r(u)|$.
	
	On the other hand, let $D'$ be a minimum $(s',t')$-dominating set of $T'$.
	Since $t'(w)=R$, we have $w\in D'$. Let $D=D'\cup (c(u)\cap r(u))$. 
	For the edge $uw$, we know that $|(N_{T'}[u]\cup N_{T'}[w])\cap D'|\geq s'(uw)$. This implies that $|(N_{T}[u]\cup N_{T}[w])\cap D|\geq s'(uw)+|c(u)\cap r(u)|\geq s(uw)$. As $w\in D$, for every other edge $uz$ incident on $u$ in $T$, we have $|(N_{T}[u]\cup N_{T}[z])\cap D|\geq |c(u)\cap r(u)|+1=s(uv)\geq s(uz)$. Since every other label remains the same, for every other edge $xy$, $|(N_{T}[x]\cup N_{T}[y])\cap D|\geq s(xy)$. Thus, $D$ is an $(s,t)$-dominating set in $T$. Hence, $\gamma_{(s,t)}(T)\leq\gamma_{(s',t')}(T')+|c(u)\cap r(u)|$. Therefore, we have $\gamma_{(s,t)}(T)=\gamma_{(s',t')}(T')+|c(u)\cap r(u)|$.
\end{proof}

\begin{lem}
	Let $s(uv) < |N_{T}[u]|$, $s(uv)>|r(u)|$ and $(s(uv)- |c(u)\cap r(u)|)\geq 2$. If $t(u)=B$, then $\gamma_{(s,t)}(T)=\gamma_{(s',t')}(T')+s(uv)-2$, where $T'$ is obtained from $T$ by deleting $c(u)$ and by relabelling $t'(u)=t'(w)=R$ and $s'(uw)= \max\{(s(uw)-s(uv)+2),0\}$ and every other label remains the same.
\end{lem}
\begin{proof}
	Let $D$ be a minimum $(s,t)$-dominating set of $T$. Since $s(uv)-|c(u)\cap r(u)|\geq 2$, at least $2$ vertices of $N_{T}[u]\setminus (c(u)\cap r(u))$, say $x$ and $y$, are in $D$. Suppose that $\{u,w\}\in D$. Let $D'=D\setminus c(u)$. We show that $D'$ is an $(s',t')$-dominating set of $T'$. Clearly, $u,w\in D'$. For the edge $uw$, we know that $|(N_{T}[u]\cup N_{T}[w])\cap D|\geq s(uw)$. This implies that $|(N_{T'}[u]\cup N_{T'}[w])\cap D'|\geq \max\{(s(uw)-s(uv)+2),0\}=s'(uw)$. Since every other label remains the same, for every other edge $xy$, we have $|(N_{T'}[x]\cup N_{T'}[y])\cap D'|\geq s'(xy)$. Hence, $D'$ is an $(s',t')$-dominating set of $T'$. Therefore, $\gamma_{(s',t')}(T')\leq \gamma_{(s,t)}(T)-|c(u)\cap r(u)|$. Now, let us assume that $\{u,w\}\notin D$. Let $D'=(D\setminus \{x,y\})\cup \{u,w\}$. Clearly, $D'$ is also an $(s,t)$-dominating set of $T$. Let $D''=D'\setminus c(u)$. We show that $D''$ is an $(s',t')$-dominating set of $T'$. Clearly, $u,w\in D''$. For the edge $uw$, we know that $|(N_{T}[u]\cup N_{T}[w])\cap D'|\geq s(uw)$. This implies that $|(N_{T'}[u]\cup N_{T'}[w])\cap D''|\geq \max\{(s(uw)-s(uv)+2),0\}=s'(uw)$. Since every other label remains the same, for every other edge $xy$, we have $|(N_{T'}[x]\cup N_{T'}[y])\cap D''|\geq s'(xy)$. Hence, $D''$ is an $(s',t')$-dominating set of $T'$. Therefore, $\gamma_{(s',t')}(T')\leq \gamma_{(s,t)}(T)-|c(u)\cap r(u)|$.
	
	Let $D'$ be a minimum $(s',t')$-dominating set of $T'$.
	Since $t'(u)=t'(w)=R$ in $T'$, we have $u,w\in D'$. Let $D=(D'\cup (c(u)\cap r(u)))\cup l(u)$, where $l(u)$ is a subset of $(c(u)\setminus r(u))$ of size $s(uv)-|c(u)\cap r(u)|-2$. Clearly, $u,w\in D$. For the edge $uw$, we know that $|(N_{T'}[u]\cup N_{T'}[w])\cap D'|\geq s'(uw)$. This implies that $|(N_{T}[u]\cup N_{T}[w])\cap D|\geq s'(uw)+s(uv)-2\geq s(uw)$. Also, for every other edge $uz$ incident on $u$, we have $|(N_{T}\cup N_{T}[z])\cap D|\geq s(uv)\geq s(uz)$. Since every other label remains the same, for every other edge $xy$, we have $|(N_{T}[x]\cup N_{T}[y])\cap D|\geq s(xy)$. Therefore, $D$ is an $(s,t)$-dominating set of $T$. Hence, $\gamma_{(s,t)}(T)\leq \gamma_{(s',t')}(T')+s(uv)-2$. Thus, $\gamma_{(s,t)}(T)=\gamma_{(s',t')}(T')+s(uv)-2$.
\end{proof}

\begin{lem}
	Let $s(uv) < |N_{T}[u]|$ and $s(uv)>|r(u)|$. If $t(u)=R$ and $(s(uv)-|c(u)\cap r(u)|)\geq 2$, then $\gamma_{(s,t)}(T)=\gamma_{(s',t')}(T')+s(uv)-2$ where $T'$ is obtained from $T$ by deleting $c(u)$ and by relabelling $t'(w)=R$ and $s'(uw)= max\{(s(uw)-s(uv)+2),0\}$ .
\end{lem}
\begin{proof}
	Let $D$ be a minimum $(s,t)$-dominating set of $T$. Since $s(uv)-|c(u)\cap r(u)|\geq 2$ and $t(u)=R$, at least $2$ vertices of $N_{T}[u]\setminus (c(u)\cap r(u))$, say $x$ and $y$, is contained in $D$ and $u\in D$. Without loss of generality let us assume that $x\in c(u)\setminus r(u)$. Let $D'=(D\setminus \{x\})\cup \{w\}$. Clearly, $D'$ is a minimum $(s,t)$-dominating set of $T$. Let $D''=D'\setminus c(u)$. We show that $D''$ is an $(s',t')$-dominating set of $T'$. Clearly, $u,w\in D''$. For the edge $uw$, we know that $|(N_{T}[u]\cup N_{T}[w])\cap D'|\geq s(uw)$. This implies that $|(N_{T'}[u]\cup N_{T'}[w])\cap D''|\geq \max\{(s(uw)-s(uv)+2),0\}=s'(uw)$. Since every other label remains the same, for every other edge $xy$, we have $|(N_{T'}[x]\cup N_{T'}[y])\cap D''|\geq s'(xy)$. Hence, $D''$ is an $(s',t')$-dominating set of $T'$. Therefore, $\gamma_{(s',t')}(T')\leq \gamma_{(s,t)}(T)-|c(u)\cap r(u)|$.
	
	Let $D'$ be a minimum $(s',t')$-dominating set of $T'$.
	Since $t'(w)=R$ and $t(u)=t'(u)=R$ in $T'$, we have $u,w\in D'$. Let $D=(D'\cup (c(u)\cap r(u)))\cup l(u)$ where $l(u)$ is a subset of $c(u)\setminus r(u)$ and $|l(u)|=s(uv)-|c(u)\cap r(u)|-2$. Clearly, $u,w\in D$. For the edge $uw$, we know that $|(N_{T'}[u]\cup N_{T'}[w])\cap D'|\geq s'(uw)$. This implies that $|(N_{T}[u]\cup N_{T}[w])\cap D|\geq s'(uw)+s(uv)-2\geq s(uw)$. Also, for every other edge $uz$ incident on $u$, we have $|(N_{T}[u]\cup N_{T}[z])\cap D|\geq s(uv)\geq s(uz)$. Since every other label remains the same, for every other edge $xy$, we have $|(N_{T}[x]\cup N_{T}[y])\cap D|\geq s(xy)$. Therefore, $D$ is an $(s,t)$-dominating set of $T$. Hence, $\gamma_{(s,t)}(T)\leq \gamma_{(s',t')}(T')+s(uv)-2$. Thus, $\gamma_{(s,t)}(T)=\gamma_{(s',t')}(T')+s(uv)-2$.
\end{proof}


All the above lemmas show that Algorithm \ref{Algo:k-VEDS_tree} returns a minimum $(s,t)$-dominating set of a given tree $T$. Now we analyze the running time of Algorithm \ref{Algo:k-VEDS_tree}. The vertex ordering in line $1$ can be computed in $O(n)$ time. For a support vertex $u$, we can find the child of $u$, say $v$, such that $uv$ has maximum $s$-label in $O(deg(u))$ time. The set of neighbours of $u$ having $t$-label as $R$ can also be computed in $O(deg(u))$ time. In each cases within the for loop in line $3 - 35$, we update the so far constructed $(s,t)$-dominating set $D$ and update the $s$-label and $t$-labels of constant number of edges and vertices from the neighbourhood of $u$. This also takes $O(deg(u))$ time. Therefore, Algorithm \ref{Algo:k-VEDS_tree} takes linear time to execute as sum of degrees is linear. As mentioned earlier, if $t(v)=B$ for every $v\in V$ and $s(e)=k$ for every $e\in E$, then Algorithm \ref{Algo:k-VEDS_tree} output a minimum $k$-vertex edge dominating set. Therefore, we have the following main theorem of this section:
%
%

\begin{theo}
	The \textsc{Min$k$VEDP} can be solved in linear time for trees.
\end{theo}
\section{Approximation algorithm and hardness}\label{sec4}

\subsection{Upper bound on approximation ratio}\label{subsec4.1}
In this subsection, we describe an approximation algorithm for \textsc{Min$k$VEDP}. This approximation algorithm follows from the existing approximation algorithm of a generalization of the classical set cover problem. The general set cover problem is defined as follows:

\noindent\underline{\textbf{General Set Cover} (\textsc{GenSetCover})}

\noindent\emph{Instance}: A set $X$, a family $\mathcal{F}$ of subsets of $X$ and an integer $k$.

\noindent\emph{Solution}: A $k$-cover of $X$, that is, a subfamily $\mathcal{C}$ of $\mathcal{F}$ such that for every $x\in X$, there are at least $k$ sets in $C$ containing $x$.

\noindent\emph{Measure}: Cardinality of the $k$-cover $|\mathcal{C}|$.


In \cite{KLASING200475}, the authors proposed an approximation algorithm for solving \textsc{GenSetCover} problem whose approximation ratio is $\ln(|F_m|)+1$, where $F_m$ is a set in $\mathcal{F}$ of maximum cardinality. Our goal is to reduce an instance of \textsc{Min$k$VEDP} into an instance of \textsc{GenSetCover} and apply the approximation algorithm for this new instance. The reduction is as follows: given a graph $G=(V, E)$, we take $X=E$. For every $v\in V$, we define the set $F_v=\{e\in E: e \ is\  incident\  to\  a~ vertex ~ in ~N_{G}[v]\}$. We set $\mathcal{F}=\{F_v|v\in V\}$. Let the approximation algorithm proposed in \cite{KLASING200475} returns a $k$-cover $\mathcal{C}$ and $D=\{v\in V| F_v\in \mathcal{C}\}$. It is easy to observe that $D$ is a $k$-ve dominating set of $G$ as $\mathcal{C}$ is a $k$-cover of $X$. Also, the maximum cardinality of a set in $\mathcal{F}$ is at most $\Delta^2(G)$, where $\Delta(G)$ is the maximum degree in $G$. Therefore, we have the following theorem:

\begin{theo}\label{Theo: Approximation algo}
	Given a graph $G=(V,E)$, the \textsc{Min$k$VEDP} can be approximated within a factor of $O(\log(\Delta(G)))$, where $\Delta(G)$ is the maximum degree in $G$.
\end{theo} 


\subsection{Lower bound on approximation ratio}\label{subsec4.2}

In this subsection, we prove a lower bound on the approximation ratio for \textsc{Min$k$VEDP} by reducing an instance of minimum vertex-edge domination problem into an instance of \textsc{Min$k$VEDP}. The minimum vertex-edge domination problem is defined as follows:

\noindent\underline{\textbf{Minimum Vertex-Edge Domination Problem}(\textsc{MinVEDP})}

\noindent\emph{Instance}: A graph $G=(V,E)$.

\noindent\emph{Solution}: A vertex-edge dominating set $D$ of $G$.

\noindent\emph{Measure}: Cardinality of the vertex-edge dominating set.

Lewis \cite{lewis} proved that for a graph $G=(V,E)$, \textsc{MinVEDP} cannot be approximated within a factor of $(1-\epsilon)\ln|V|$ for any $\epsilon >0$, unless $NP\subseteq DTIME(|V|^{O(\log\log|V|)})$. Next, we show an approximation preserving reduction from \textsc{MinVEDP} to \textsc{Min$k$VEDP}. Let $G=(V,E)$ be an instance of \textsc{MinVEDP}. The construction of $G'=(V', E')$, an instance of \textsc{Min$k$VEDP} is as follows: consider a clique $C=\{c_1, c_2, \ldots, c_{k-1}\}$ on $k-1$ vertices. For every $v_i\in V$, add an edge between $v_i$ and $c_j$ for every $c_j\in C$. Consider a vertex $v$ and make $v$ adjacent to every $c_j\in C$. Finally, add another vertex $u$ which is adjacent to $v$. Therefore, $V^{'}=V\cup C\cup \{v,u\}$ and $E^{'}=E\cup \{c_ic_j|c_i,c_j\in C\}\cup\{c_iv|c_i\in C\}\cup\{vu\}$. The construction of $G'$ is illustrated in Figure \ref{Fig: NPC for chordal}. 
\begin{figure}
	\centering
	\captionsetup{justification=centering}
	\includegraphics[scale=0.65]{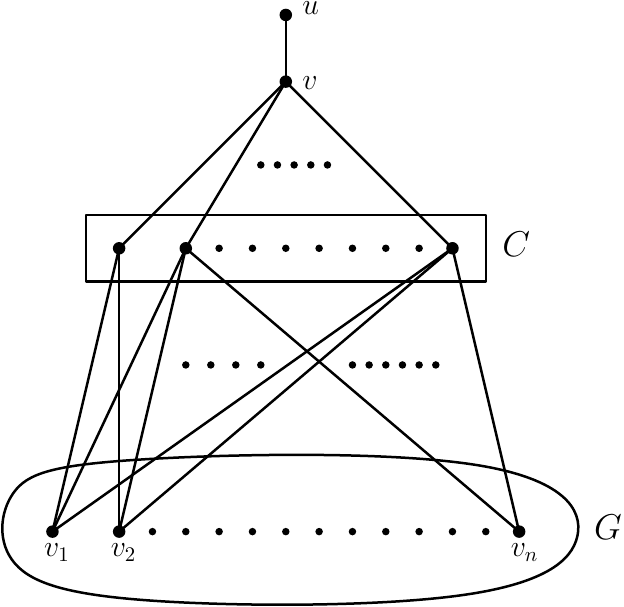}
	\caption{Construction of $G'$ from $G$}\label{Fig: NPC for chordal}
\end{figure}

\begin{cl}\label{Claim:ve-domtok-vedom}
	The graph $G$ has a ve-dominating set of size at most $t$ if and only if $G'$ has a $k$-ve dominating set of size at most $t+k$.
\end{cl}
\begin{proof}
	Let $D$ be a ve-dominating set of $G$ of size at most $t$. Consider the set $D'=D\cup C\cup \{v\}$. Note that, every edge incident to $C\cup\{v\}$ is $k$-ve dominated by $C\cup \{v\}$ and hence by $D'$. Every other edge is of the form $v_iv_j$. For such an edge $v_iv_j$, we have $C\subseteq N_{G'}[v_i]\cup N_{G'}[v_j]$ by construction and $|(N_{G}[v_i]\cup N_{G}[v_j])\cap D|\geq 1$ as $D$ is a ve-dominating set of $G$. Hence, $|(N_{G'}[v_i]\cup N_{G'}[v_j])\cap D_{k}|\geq k$. Therefore, $D'$ is a $k$-ve dominating set of $G^{'}$ of size at most $t+k$.
	
	
	On the other hand, let $D'$ be a $k$-ve dominating set of $G^{'}$ of size at most $t+k$. Since $(N_{G'}[u]\cup N_{G'}[v])= C\cup \{v,u\}$, to $k$-ve dominate $vu$ at least $k$ vertices from $C\cup \{v,u\}$ must be present in $D'$. If $C\cup \{v,u\}\subset D'$, then we can remove $u$ from $D'$ to construct another $k$-ve dominating set of $G'$ of size at most $t+k$. Further, if $u\in D'$ but $C\cup \{v\} \not\subset D'$, then we can replace $u$ by the missing vertex from $C\cup \{v\}$ to construct another $k$-ve dominating set of $G'$ of size at most $t+k$. Therefore, without loss of generality, we can assume that $D'$ is a $k$-ve dominating set of $G^{'}$ of size at most $t+k$ not containing $u$ but $C\cup \{v\} \subset D'$. Let $D= D'\setminus (C\cup \{v\})$. Clearly, $|D|$ is at most $t$. For any edge $v_iv_j\in E$, we know that $|(N_{G'}[v_i]\cup N_{G'}[v_j])\cap D'|\geq k$. Since $C\subset (N_{G'}[v_i]\cup N_{G'}[v_j]) \cap D'$, we have $|(N_{G}[v_i]\cup N_{G}[v_j])\cap D|\geq 1$. Therefore, $D$ is a ve-dominating set of $G$ of size at most $t$.
\end{proof}

%
%

Next, by using the above construction, we show the lower bound on the approximation ratio. 

\begin{theo}
	For a graph $G=(V,E)$, \textsc{Min$k$VEDP} cannot be approximated within a factor of $(1-\epsilon)\ln|V|$ for any $\epsilon >0$, unless $NP\subseteq DTIME(|V|^{O(\log\log|V|)})$. 	
\end{theo} 
\begin{proof}
	Let $\mathcal{A}$ be an approximation algorithm to find a $k$-ve dominating set whose approximation ratio is $\rho$. First note that if $k$ is a constant, then, given a graph $G$ whose minimum ve-dominating set size is at most  $k$, we can solve \textsc{MinVEDP} in polynomial time. Let us consider a graph $G$ whose minimum ve-dominating set size is more than $k$. We find an approximate ve-dominating set of $G$ as follows: first, we construct $G'$ using the construction in Claim \ref{Claim:ve-domtok-vedom}. Then using the algorithm $\mathcal{A}$, we find an approximate $k$-ve dominating set, say $D'$, of $G'$. And finally, we get an approximate ve-dominating set, say $D$, of $G$ using Claim \ref{Claim:ve-domtok-vedom}. Let $D^*$ and $D'^*$ be the minimum ve-dominating set of $G$ and minimum $k$-ve dominating set of $G'$, respectively. Therefore, we have 
	\begin{align*}
		|D| &\leq |D'|\\
		&\leq \rho |D'^*| \tag*{$[Since~ |D'|\leq \rho |D'^*|]$}\\
		&\leq \rho (|D^*|+k) \tag*{$[Since~ |D'^*|\leq |D^*|+k~ by~Claim~\ref{Claim:ve-domtok-vedom}]$}\\
		&\leq \rho (1+\frac{k}{|D^*|})|D^*|	
	\end{align*}	
	If possible, let \textsc{Min$k$VEDP} can be approximated within a factor of $(1-\epsilon)\ln|V'|$ for any $\epsilon >0$, that is, $\rho= (1-\epsilon)\ln|V'|$. Since, $G$ is a graph such that $|D^*|>k$, we set $\epsilon$ such that $\frac{k}{|D^*|} < \epsilon < 1$ and $\epsilon > \frac{1}{\sqrt{2}}$. Therefore, we have
	\begin{align*}
		|D| &\leq (1-\epsilon)\ln |V'|(1+\epsilon)|D^*|\\
		&\leq (2-2\epsilon^2)\ln |V||D^*| \tag*{$[Since~ |V'|\leq |V|^2]$}
	\end{align*}
	Since $\epsilon > \frac{1}{\sqrt{2}}$, we have $\epsilon'=2\epsilon^2-1$ is a non-zero quantity which is also less than $1$. Hence we have, $|D|\leq (1-\epsilon')\ln |V||D^*|$, which is a contradiction. Therefore, \textsc{Min$k$VEDP} cannot be approximated within a factor of $(1-\epsilon)\ln|V|$ for any $\epsilon>0$ unless $NP\subseteq DTIME(|V|^{O(\log\log|V|)})$.
\end{proof}

\subsection{APX-complete for bounded degree graphs}\label{subsec4.3}
In this subsection, we show that \textsc{Min$k$VEDP} is APX-complete for graphs with maximum degree $k+3$. We denote the \textsc{Min$k$VEDP} restricted to graphs with maximum degree $k+3$ by \textsc{Min$k$VEDP$(k+3)$}. First, we define the notion of $L$-reduction \cite{L-reduction}. Given two NP-complete optimization problem $\pi_1$ and $\pi_2$ and a polynomial time transformation $f$ from the instances of $\pi_1$ to the instances of $\pi_2$, we say $f$ is an $L$-reduction if there are positive constants $\alpha$ and $\beta$ such that for every instance $x$ of $\pi_1$,
\begin{enumerate}
	\item $opt_{\pi_2}(f(x))\leq \alpha opt_{\pi_1}(x)$;
	\item for every feasible solution $y$ of $f(x)$ with objective value $m_{\pi_2}(f(x),y)=c_2$, we can find a solution $y'$ of $x$ in polynomial time with $m_{\pi_1}(x,y')=c_1$ such that $|opt_{\pi_1}(x)-c_1|\leq \beta |opt_{\pi_2}(f(x))-c_2|$.
\end{enumerate}

To prove APX-completeness, we show an $L$-reduction from minimum $k$-tuple domination problem for graphs with maximum degree bounded by $k+2$, which is known to be APX-complete \cite{KLASING200475}, to our problem \textsc{Min$k$VEDP$(k+3)$}. The problem is defined as follows:

\noindent\underline{\textbf{Minimum $k$-tuple Domination Problem }(\textsc{Min$k$DOM$(k+2)$})}

\noindent\emph{Instance}: A graph $G=(V,E)$ of degree bounded by $(k+2)$ and $k\geq 2$.

\noindent\emph{Solution}: A $k$-tuple dominating set of $G$.

\noindent\emph{Measure}: Cardinality of the $k$-tuple dominating set.

\begin{theo}
	\textsc{Min$k$VEDP$(k+3)$} is APX-complete for $k\geq 2$.
\end{theo}
\begin{proof}
	By Theorem \ref{Theo: Approximation algo} we have if the maximum degree is bounded by a constant, then the approximation ratio is also constant. Therefore, the problem is in APX. From a graph $G=(V,E)$ with maximum degree bounded by $(k+2)$, we construct a graph $G^{'}=(V^{'}, E^{'})$, where $V^{'}=V\cup\{u_i|1\leq i\leq |V|\}$ and $E^{'}=E\cup\{v_iu_i| 1\leq i\leq |V|\}$. Clearly, the maximum degree of $G^{'}$ is bounded by $(k+3)$. 
	
	Let $D^{*}$ be a minimum $k$-tuple dominating set of $G$. For every $v_i\in V$, we have $|N_{G}[v_i]\cap D^{*}|\geq k$. Therefore, for every $v_iv_j\in E$, we have $|(N_{G}[v_i]\cup N_{G}[v_j])\cap D^{*}|\geq k$. Also, for every edge $u_iv_i$, we have $|(N_{G^{'}}[u_i]\cup N_{G^{'}}[v_i])\cap D^{*}|\geq k$. Hence, $D^{*}$ is a $k$-ve dominating set of $G^{'}$. Hence, $|D'^{*}|\leq |D^{*}|$, where $D'^{*}$ is a minimum $k$-ve dominating set of $G^{'}$. On the other hand, let $D'^{*}$ be a minimum $k$-ve dominating set of $G^{'}$. Therefore, for every edge $v_iu_i$, we have $|(N_{G^{'}}[u_i]\cup N_{G^{'}}[v_i])\cap D'^{*}|\geq k$. If $|N_{G}[v_i]\cap D'^{*}|\geq k$ for every vertex $v_i\in V$, then $D'^{*}$ is a $k$-tuple dominating set. Otherwise, if $u_i\in D'^{*}$, then $D=(D'^{*}\setminus \{u_i\})\cup\{v_z\}$ is a $k$-tuple dominating set where $v_z\in N_{G}[v_i]\setminus D'^{*}$. Therefore, $|D^*|\leq |D'^*|$, where $D^*$ is a minimum $k$-tuple dominating set of $G$. Hence, we have $|D^{*}|=|D^{'*}|$.

	Similarly, we can show that $D$ is $k$-tuple dominating set of $G$ if and only if $D$ is $k$-ve dominating set of $G'$. Therefore, the above reduction is an $L$-reduction with $\alpha=1$ and $\beta=1$. Hence, \textsc{Min$k$VEDP$(k+3)$} is APX-complete.
\end{proof}

\section{Conclusion}\label{sec5}
In this paper, we have shown that the \textsc{Decide$k$VEDP} is NP-complete for chordal graphs. In algorithmic point of view \textsc{Min$k$VEDP}, we have designed a linear time algorithm in tree. Also, we introduce an approximation algorithm for \textsc{Min$k$VEDP}, establish  the lower bound of approximation ratio for the same and showed that \textsc{Min$k$VEDP} is APX-complete in bounded degree graphs. 
It would be interesting to study the complexity status of this problem in different subclasses of chordal graphs. 



\section*{Declarations}

\begin{itemize}
	\item Funding: Subhabrata Paul is supported by CSIR-HRDG Research Grant (No. 25(0313)/20/EMR-II).
	
	\item Competing interests: The authors declare that they have no conflict of interest.
	
	\item Availability of data and materials: This work has no associated data.
	
	\item Authors' contributions: All authors contributed equally to this work.
\end{itemize}

\bibliographystyle{alpha}
\bibliography{VEDom_bib}

\begin{thebibliography}{NKPV21}

\bibitem[BCHH16]{boutrig2016vertex}
R.~Boutrig, M.~Chellali, T.~W Haynes, and S.T. Hedetniemi.
\newblock Vertex-edge domination in graphs.
\newblock {\em Aequationes mathematicae}, 90:355--366, 2016.

\bibitem[Cha98]{labelingC}
Gerard~J. Chang.
\newblock {\em Algorithmic Aspects of Domination in Graphs}, pages 1811--1877.
\newblock Springer US, Boston, MA, 1998.

\bibitem[CS12]{chitra2012global}
S.~Chitra and R.~Sattanathan.
\newblock Global vertex-edge domination sets in graph.
\newblock In {\em Int. Math. Forum}, volume~7, pages 233--240, 2012.

\bibitem[CS22]{chen2022double}
Xue-G. Chen and M.Y. Sohn.
\newblock Double vertex-edge domination in trees.
\newblock {\em Bulletin of the Korean Mathematical Society}, 59(1):167--177,
  2022.

\bibitem[GJ79]{garey1979}
M.R. Garey and D.S. Johnson.
\newblock {\em Computers and intractability:A Guide to the Theory of
  NP-Completeness}.
\newblock W.H Freeman, New York, 1979.

\bibitem[HHS98a]{domingraphbook}
T.W. Haynes, S.~Hedetniemi, and P.~Slater.
\newblock {\em Domination in Graphs}, volume~2.
\newblock New York, 1998.

\bibitem[HHS98b]{Fundamentalbook}
T.W. Haynes, S.~Hedetniemi, and P.~Slater.
\newblock {\em Fundamentals of Domination in Graphs (1st ed.)}.
\newblock New York, 1998.

\bibitem[JD22]{jena}
S.K. Jena and G.K. Das.
\newblock Vertex-edge domination in unit disk graphs.
\newblock {\em Discrete Applied Mathematics}, 319:351--361, 2022.

\bibitem[KL04]{KLASING200475}
Ralf Klasing and Christian Laforest.
\newblock Hardness results and approximation algorithms of k-tuple domination
  in graphs.
\newblock {\em Information Processing Letters}, 89(2):75--83, 2004.

\bibitem[KVK14]{krishnakumari2014bounds}
B.~Krishnakumari, Y.B. Venkatakrishnan, and M.~Krzywkowski.
\newblock Bounds on the vertex--edge domination number of a tree.
\newblock {\em Comptes rendus mathematique}, 352(5):363--366, 2014.

\bibitem[Lew07]{lewis}
J.R. Lewis.
\newblock {\em Vertex-edge and edge-vertex parameters in graphs.}
\newblock PhD thesis, Clemson University, Clemson, SC, USA, 2007.

\bibitem[LW23]{li2023polynomial}
Peng Li and Aifa Wang.
\newblock Polynomial time algorithm for k-vertex-edge dominating problem in
  interval graphs.
\newblock {\em Journal of Combinatorial Optimization}, 45(1):45, 2023.

\bibitem[NKPV21]{naresh}
H.~Naresh~Kumar, D.~Pradhan, and Y.B. Venkatakrishnan.
\newblock Double vertex-edge domination in graphs: complexity and algorithms.
\newblock {\em Journal of Applied Mathematics and Computing}, 66(1):245--262,
  2021.

\bibitem[Pet86]{peters}
K.W.J. Peters.
\newblock {\em Theoritical and algorithmic results on domination and
  connectivity}.
\newblock PhD thesis, Clemson University, Clemson, SC, USA, 1986.

\bibitem[PPV21]{paul2}
S.~Paul, D.~Pradhan, and S.~Verma.
\newblock Vertex-edge domination in interval and bipartite permutation graphs.
\newblock {\em Discussiones Mathematicae: Graph Theory}, 2021.

\bibitem[PR21]{paul}
S.~Paul and K.~Ranjan.
\newblock Results on vertex-edge and independent vertex-edge domination.
\newblock {\em Journal of Combinatorial Optimization}, pages 1--28, 2021.

\bibitem[PY88]{L-reduction}
Christos Papadimitriou and Mihalis Yannakakis.
\newblock Optimization, approximation, and complexity classes.
\newblock In {\em Proceedings of the Twentieth Annual ACM Symposium on Theory
  of Computing}, STOC '88, page 229–234, New York, NY, USA, 1988. Association
  for Computing Machinery.

\bibitem[{\.Z}yl19]{zylinski2019vertex}
P.~{\.Z}yli{\'n}ski.
\newblock Vertex-edge domination in graphs.
\newblock {\em Aequationes mathematicae}, 93(4):735--742, 2019.

\end{thebibliography}

\end{document}